\documentclass[a4paper]{amsart}

\newtheorem{theorem}{Theorem}[section]
\newtheorem{lemma}[theorem]{Lemma}

\newtheorem{definition}[theorem]{Definition}
\newtheorem{example}[theorem]{Example}

\newcommand{\bin}[2]{\genfrac{(}{)}{0pt}{}{#1}{#2}}
\newcommand{\on}[2]{\genfrac{}{}{0pt}{}{#1}{#2}}
\newcommand{\R}{\mathbb{R}}
\newcommand{\C}{\mathbb{C}}
\newcommand{\Sn}{\mathfrak{S}_n}
\newcommand{\n}{\mathfrak{n}}
\newcommand{\E}{\mathbb{E}}
\newcommand{\Cln}{\mathrm{Cl}_n}
\newcommand{\type}{\mathrm{type}}
\newcommand{\stat}{\mathfrak{s}}
\newcommand{\cyc}{\mathrm{cyc}}
\newcommand{\exc}{\mathrm{exc}}
\newcommand{\des}{\mathrm{des}}
\newcommand{\maj}{\mathrm{maj}}
\newcommand{\id}{\mathrm{id}}
\newcommand{\inv}{\mathrm{inv}}
\newcommand{\wexc}{\mathrm{wexc}}
\newcommand{\Ht}{\mathrm{ht}}
\newcommand{\av}[1]{\overline{#1}}

\begin{document}
\author{Axel Hultman}
\address{Department of Mathematics, Link\"oping University, SE-581 83, Link\"oping, Sweden.}
\title[Statistics of random products]{Permutation statistics of products of random
  permutations}

\begin{abstract}
Given a permutation statistic $\stat : \Sn \to \R$, define the {\em
  mean statistic} $\av{\stat}$ as the statistic which computes the
mean of $\stat$ over conjugacy classes. We describe a way to calculate
the expected value of $\stat$ on a product of $t$ independently chosen
elements from the uniform distribution on a union of conjugacy classes
$\Gamma\subseteq \Sn$. In order to apply the formula, one needs to
express the class function $\av{\stat}$ as a linear combination of irreducible
$\Sn$-characters. We provide such expressions for several commonly
studied permutation statistics, including the excedance number, inversion
number, descent number, major index and $k$-cycle number. In particular,
this leads to formulae for the expected values of said statistics.
\end{abstract}

\maketitle

\section{Introduction}
Consider the symmetric group $\Sn$ of permutations of $[n] = \{1,
\ldots, n\}$. For $\Gamma \subseteq \Sn$ one may study the behaviour of various permutation
statistics $\stat : \Sn \to \R$ on products $\gamma_1\cdots \gamma_t \in \Sn$ of
random $\gamma_i \in \Gamma$.
\begin{definition}
Choose a subset $\Gamma \subseteq \Sn$, a function $\stat : \Sn \to
\R$ and a nonnegative integer $t$. We denote by $\E_\Gamma(\stat, t)$
the expected value of $\stat$ on a product of $t$ elements
independently chosen from the uniform distribution on $\Gamma$.
\end{definition}

A product of $t$ random elements of $\Gamma$ corresponds to a $t$-step random
walk on the Cayley
graph of $\Sn$ induced by $\Gamma$. Random walks on Cayley graphs form
a classical and well studied subject in probability
theory; a good general reference is \cite{SC}. In the present paper, we specifically
address the problem of computing $\E_\Gamma(\stat, t)$. Recent work in
this vein includes the following. The case of $\Gamma$ being the set of adjacent
transpositions and $\stat$ counting inversions was studied by Eriksson
et al.\ \cite{EES}, Eriksen \cite{eriksen} and Bousquet-M\'elou
\cite{BM}. Turning, instead, to $\Gamma$ comprised of all
transpositions, $\stat$ being the absolute length function
(essentially counting disjoint cycles) was studied in \cite{EH}, whereas
J\"onsson \cite{jonsson} considered the fixed point counting function
$\stat$ and Sj\"ostrand \cite{sjostrand} found the solution when
$\stat$ counts inversions.  

We shall describe a method to attack the general
problem. Although it could potentially be of use for more general
$\Gamma$ (as indicated by the hypotheses of Theorem \ref{th:inner}
below), we shall apply it to situations where $\Gamma$ is a union of 
conjugacy classes. The technique makes use of representation
theory of $\Sn$. Similar ideas have been frequent in the study
of random walks on Cayley graphs since the seminal paper by Diaconis and Shahshahani \cite{DS}. In particular, our
method is heavily inspired by that described in
\cite{EH}. We are concerned with more general $\Gamma$, but the principal novelty here is to dispose of the apparent
requisite of \cite{EH} that $\stat$ be a class function. As we shall
see, removing this restriction significantly improves the versatility of the method. 

Here is a brief sketch of the content of the paper. In the next
section, we review some facts about symmetric group characters. In
Section \ref{se:main}, we then describe how they connect with our
expected statistics problem. The main results of that section, Theorem \ref{th:inner} and
Theorem \ref{th:blambda}, provide a recipe for computing
$\E_\Gamma(\stat, t)$ whenever $\Gamma$ is a union of conjugacy
classes. In order to express $\E_\Gamma(\stat, t)$ explicitly for a given
statistic $\stat$, the remaining task is to decompose the {\em mean statistic} $\av{\stat}$ as a
linear combination of irreducible $\Sn$-characters. This turns out to be a fairly straightforward task for many
standard permutation statistics. We provide explicit decompositions
for the mean statistics corresponding to the $k$-cycle number, excedance number, weak
excedance number, inversion number, major index and descent number
statistics in Sections \ref{se:cycle}, \ref{se:exc} and
\ref{se:inv}. In particular, $\E_\Gamma(\stat, t)$ is determined for $\stat$ being
any of these statistics (and conjugation invariant $\Gamma$). We conclude
with a few explicit examples in Section \ref{se:examples}. 
  
\section{Symmetric group characters}\label{se:character}
In this section, we review elements of the representation theory of
$\Sn$. From this vast and classical subject only a few
bits and pieces that we need in the sequel are extracted in order to agree on notation. We
refer to \cite{sagan} for a thorough background and much more information.

Let $P_n$ denote the set of integer partitions $\lambda \vdash n$. The
irreducible representations of $\Sn$ are in bijection with $P_n$ in a
standard way. We use
the notation $\rho^\lambda$ for the representation indexed by $\lambda
\vdash n$ and denote the corresponding character by
$\chi^\lambda$. These irreducible characters form a basis for the
$\C$-vector space $\Cln = \{f:P_n \to \C\}$ of class functions. Moreover, this basis is
orthonormal with respect to the standard Hermitian inner product on $\Cln$
\[
\langle f,g\rangle = \frac{1}{n!}\sum_{\lambda \vdash n} |C_\lambda |f(\lambda)g(\lambda)^*,
\] 
where $C_\lambda$ denotes the conjugacy class of permutations $\pi$ with
$\type (\pi) = \lambda$ as cycle type. 

Abusing notation, we at times consider class functions as defined on
$\Sn$ rather than on $P_n$. In other words, for $\pi \in \Sn$ and $f\in \Cln$,
$f(\pi)$ should be interpreted as $f(\type (\pi))$. We trust the
context to prevent confusion.

It is convenient to encode partitions as weakly
decreasing sequences of positive integers, sometimes employing exponent notation
to signal repeated parts. For example, $(7, 3^4, 1)$ denotes the
partition of $20$ which consists of one part of size $7$, four parts of
size $3$ and one part of size $1$. In this notation, a {\em hook shape} is a partition
of the form $(a, 1^b)$ for integers $a\ge 1$ and $b\ge 0$.

The trivial $\Sn$-character is indexed by $(n)$. Thus, $\chi^{(n)}(\mu) = 1$
for all $\mu \vdash n$. The next three lemmata collect a few more values of certain irreducible
characters that we shall need in the
sequel. All statements are readily verified using e.g.\ the Murnaghan-Nakayama rule. 

Define $f^\lambda = \chi^\lambda((1^n))$. This
is the dimension of the irreducible representation $\rho^\lambda$.

\begin{lemma}\label{le:hook}
Hook shape characters satisfy
\[
f^{(n-k, 1^k)} = \bin{n-1}{k}
\]
and
\[
\chi^{(n-k,1^k)}((n)) = (-1)^k,
\]
whereas $\chi^\lambda((n)) = 0$ if $\lambda$ is not a hook shape.
\end{lemma}

\begin{lemma}\label{le:explicit}
Let $\lambda \vdash n$. If $\lambda$ has $p$ parts of size $1$ and $q$
parts of size $2$, then
\[
\chi^{(n-1,1)} (\lambda) = p-1
\]
and
\[
\chi^{(n-2,1^2)} (\lambda) = \bin{p-1}{2}-q.
\]
\end{lemma}

A frequently occurring quantity is the {\em content} of $\lambda
\vdash n$. It is defined by
\[
c_\lambda = \bin{n}{2}\frac{\chi^\lambda((2, 1^{n-2}))}{f^\lambda}.
\]

\begin{lemma}\label{le:content}
The content of a hook shape is
\[
c_{(n-k, 1^k)} = \frac{n(n-2k-1)}{2}.
\]
\end{lemma}

\section{Mean statistics and irreducible characters}\label{se:main}
Let $\stat:\Sn \to \R$ be any real-valued function on the symmetric
group. In all our subsequent applications, $\stat$ will be a permutation
statistic associating a nonnegative integer with each permutation
in $\Sn$.

Choose $\Gamma\subseteq \Sn$. We now describe a general procedure for
computing the expected value of $\stat$ on a product of random
elements of $\Gamma$ taken independently from the uniform
distribution.

The {\em mean statistic} $\av{\stat}$ is the class function which computes the mean
of $\stat$ over conjugacy classes. That is, $\av{\stat}:P_n\to \R$ is
defined by
\[
\av{\stat}(\lambda) = \frac{1}{|C_\lambda|}\sum_{\pi\in C_\lambda} \stat(\pi). 
\]
Hence, thinking of class functions as being defined on $\Sn$, $\av{\stat}
= \stat$ if and only if $\stat$ is a class function. Considered in that particular setting, the
remainder of this section resembles the procedure described in
\cite{EH}.\footnote{The context of \cite{EH} was that of $\stat(\pi)$
  being the absolute length of $\pi$ with $\Gamma$ the  
set of transpositions. The technique, however, could clearly have been
applied to any class function $\stat$.}

Let $\pi \in \Sn$. We need to keep track of the ways to express $\pi$ as a product of elements from $\Gamma$. To this end, define a permutation statistic $\n_t$ by
\[
\n_t(\pi) = \#\{(\gamma_1, \ldots, \gamma_t) \in \Gamma^t\mid \gamma_1\cdots \gamma_t = \pi\}.
\]
Observe that all $\n_t$ are class functions if and only if $\Gamma$ is a union of conjugacy classes.

\begin{theorem}\label{th:inner}
If at least one of the statistics $\n_t$ and $\stat$ is a class function, then
\[
\E_\Gamma (\stat, t) = \frac{n!}{|\Gamma|^t}\langle  \av{\stat}, \av{\n_t} \rangle.
\] 
\begin{proof}
By definition,
\[
\E_\Gamma (\stat, t) = \frac{1}{|\Gamma|^t}\sum_{\pi\in\Sn}\n_t(\pi)\stat(\pi).
\]
Assume now that $\stat$ is a class function. Rewriting the right hand side by first summing over the conjugacy classes of $\Sn$, we obtain
\[
\begin{split}
\E_\Gamma (\stat, t) &= \frac{1}{|\Gamma|^t}\sum_{\lambda \vdash
  n}\av{\stat}(\lambda)\sum_{\pi\in C_\lambda}\n_t(\pi)\\
 &= \frac{1}{|\Gamma|^t}\sum_{\lambda \vdash n}\av{\stat}(\lambda)|C_\lambda|\av{\n_t}(\lambda)\\
 &= \frac{n!}{|\Gamma|^t}\langle \av{\stat}, \av{\n_t}\rangle,
\end{split}
\]
as desired. If, instead, $\n_t$ is a class function, the proof is completely analogous.
\end{proof}
\end{theorem}

Under the hypotheses of the preceding theorem, we are left
with the task of evaluating the inner product of two mean
statistics. This is easy if we are somehow able to express them in the
orthonormal basis comprised of the irreducible $\Sn$-characters. In other words, we want to find the coefficients $a_\lambda$ and $b_\lambda^{(t)}$ defined by
\[
\av{\stat} = \sum_{\lambda \vdash n} a_\lambda \chi^\lambda
\]
and
\[
\av{\n_t} = \sum_{\lambda \vdash n} b_\lambda^{(t)} \chi^\lambda,
\]
respectively. 

Although Theorem \ref{th:inner} applies if $\n_t$ {\em or} $\stat$ is a
class function, all our subsequent applications come from the former
setting. The next theorem is the reason; it shows how to compute
the $b_\lambda^{(t)}$ if $\Gamma$ consists of conjugacy
classes. Variations of the formula (and its proof) are numerous in the
literature. With a bit of labour, it can be extracted e.g.\ from \cite{DS} or \cite{ito}. When
$\Gamma$ is a conjugacy class, the statement follows immediately from
\cite[Theorem A.1.9]{LZ} which is attributed to Frobenius. We provide
a self-contained proof for convenience.

\begin{theorem} \label{th:blambda}
Suppose $\Gamma$ is a disjoint union of conjugacy classes $\Gamma_1,
\ldots, \Gamma_k \subseteq \Sn$. Let $\mu_i$ denote the cycle type of
the permutations in $\Gamma_i$. Then,
\[
b_\lambda^{(t)} = \frac{1}{n!(f^\lambda)^{t-1}}\left(\sum_{i=1}^k|\Gamma_i|\chi^\lambda (\mu_i)\right)^t.
\]
\begin{proof}
Let $\delta_{\cdot, \cdot}$ denote the Kronecker delta. By the Schur orthogonality relations,
\[
\frac{1}{n!}\sum_{\lambda \vdash n} f^\lambda \chi^\lambda(\mu) = \delta_{(1^n), \mu} = \av{\n_0}(\mu).
\]
This proves the $t=0$ case of the asserted statement.

Suppose $f$ is any class function and define a linear map $P_\Gamma$ on $\Cln$ by declaring
\[
P_\Gamma(f)(\pi) = \sum_{\gamma\in \Gamma} f(\pi\gamma^{-1})
\]
for $\pi \in \Sn$. By definition, $\av{\n_t} = P_\Gamma^t(\av{\n_0})$. Thus, it suffices to
show that $\chi^\lambda$ is an eigenfunction of $P_\Gamma$ with
eigenvalue $\frac{1}{f^\lambda}\sum_{i=1}^k|\Gamma_i|\chi^\lambda
(\mu_i)$. To this end, we may by linearity assume without loss of generality that
$\Gamma$ is a single conjugacy class. Denote its cycle type simply
by $\mu$.

Define $M_\lambda = \sum_{\gamma\in \Gamma} \rho^\lambda
(\gamma)$. Since $\Gamma$ is a conjugacy class, $M_\lambda$ and
$\rho^\lambda(\pi)$ commute for all $\pi\in \Sn$. By Schur's Lemma,
$M_\lambda = h(\lambda)I$, where $I$ is the identity map and
$h(\lambda)\in \C$. Now observe that
\[
P_\Gamma(\chi^\lambda)(\pi) = \mathrm{trace} \left(\rho^\lambda(\pi)
  M_\lambda\right) = h(\lambda) \chi^\lambda(\pi)
\]
for $\pi \in \Sn$. Plugging in $\pi = \id$ yields
$|\Gamma|\chi^\lambda(\mu) = h(\lambda) f^\lambda$, proving the claim.
\end{proof}
\end{theorem}




For $\Gamma$ satisfying the hypothesis of Theorem \ref{th:blambda}, we conclude
that the remaining challenge is to decompose the mean statistic $\av{\stat}$ as a linear
combination of irreducible $\Sn$-characters. The upcoming three sections
are essentially devoted to such computations.

\section{Cycle numbers}\label{se:cycle}

For $\pi\in \Sn$ and a positive integer $k$, let $\cyc_k (\pi)$ denote
the number of elements that are contained in $k$-cycles
in the disjoint cycle decomposition of $\pi$. Thus, $\pi$ contains
$\cyc_k(\pi)/k$ $k$-cycles. Clearly, $\cyc_k$ is a class function so
that $\av{\cyc_k} = \cyc_k$. The following result, providing a decomposition of this statistic, is due to Alon
and Kozma \cite{AK}. We take this opportunity to state a 
shorter, independent proof.\footnote{The proof employs standard terminology
  from the theory of symmetric functions.  We refrain from reproducing the
  definitions since they are not used elsewhere in the
  paper. Everything can be found e.g.\ in \cite[Chapter 7]{stanley}.}

\begin{theorem}[Theorem 3 in \cite{AK}] \label{th:cycles}
Let $k\in [n]$. Regarded as a class function on $\Sn$, $\cyc_k$ decomposes as 
\[
\begin{split}
\cyc_k = \chi^{(n)} &+
\sum_{i=1}^{\min(k,n-k)}(-1)^{k-i}\chi^{(n-k,i,1^{k-i})} \\
&+ \sum_{j=n-k+1}^{k-1}(-1)^{k-j}\chi^{(j,n-k+1,1^{k-j-1})}.
\end{split}
\]
\begin{proof} If $\mu\vdash n-k$ is obtained from $\lambda\vdash n$ by removing a part of
size $k$, we of course have $\cyc_k(\lambda) = \cyc_k(\mu) + k$. Thus,
the symmetric function image of $\cyc_k$ under the characteristic map 
is 
\[
\sum_{\on{\lambda\vdash n}{\lambda \text{ has $k$-parts}}}
\frac{p_\lambda\cyc_k(\lambda)}{z_\lambda} = \sum_{\mu \vdash
  n-k}\frac{p_\mu p_k(\cyc_k(\mu) + k)}{z_\mu k(\cyc_k(\mu)/k + 1)} =
p_k\sum_{\mu\vdash n-k}\frac{p_\mu}{z_\mu} = p_ks_{n-k}.
\]
Using \cite[7.72]{stanley}, we may write
\[
p_ks_{n-k} = \sum (-1)^{\Ht(\lambda/(n-k))}s_\lambda,
\]
where the sum is over all partitions $\lambda\vdash n$ such that
$\lambda/(n-k)$ is a border strip, and $\Ht(\lambda/(n-k))$ is one
less than the number of rows in the strip. Applying the inverse of the
characteristic map, this is precisely the desired result. 
\end{proof}
\end{theorem}

Theorem \ref{th:cycles} refines results
from \cite{EH} and \cite{jonsson}. The former work essentially
revolved around the total
number of cycles, i.e.\ the statistic $\sum_k \cyc_k/k$, whereas the
fixed point number $\cyc_1$ was studied in the latter. 

\section{Excedances}\label{se:exc}
Recall that an {\em excedance} of $\pi \in \Sn$ is an index $i \in [n]$ such
that $\pi(i) > i$. Similarly, $i$ is a {\em weak excedance} if $\pi(i)
\ge i$. Let $\exc(\pi)$ and $\wexc(\pi)$ denote the number of
excedances and weak excedances, respectively, of $\pi$. Clearly,
neither $\exc$ nor $\wexc$ is a class function. 
\begin{theorem}\label{th:exc}
The mean statistics $\av{\exc}$ and $\av{\wexc}$ decompose as
\[
\av{\exc} = \frac{n-1}{2}\chi^{(n)} - \frac{1}{2}\chi^{(n-1,1)}
\]
and
\[
\av{\wexc} = \frac{n+1}{2}\chi^{(n)} + \frac{1}{2}\chi^{(n-1,1)},
\]
respectively.
\begin{proof}
If $i$ is not a fixed point of $\pi$, then $i$ is an excedance of
$\pi$ if and only if $\pi(i)$ is not an excedance of $\pi^{-1}$. A fixed
point is a weak excedance but not an excedance. Hence,
\[
\av{\exc}(\lambda) = \frac{1}{2\#C_\lambda}\sum_{\pi\in
  C_\lambda}(\exc(\pi) + \exc(\pi^{-1})) = \frac{n-p}{2},
\]
where $p$ is the number of fixed points of any $\pi \in C_\lambda$,
i.e.\ the number of parts that equal one in $\lambda$. Similarly,
\[
\av{\wexc}(\lambda) = \frac{n+p}{2}.
\]
The result now follows from Lemma \ref{le:explicit}.
\end{proof}
\end{theorem}

\section{Inversions, descents and the major index}\label{se:inv}
This section treats the mean statistics associated with
three commonly occurring permutation statistics. First, we recall
their definitions.

Let $\pi \in \Sn$. A {\em descent} of $\pi$ is an index
$i\in [n-1]$ such that $\pi(i) > \pi(i+1)$. The number of descents of
$\pi$ is denoted by $\des(\pi)$, whereas the {\em major index}
$\maj(\pi)$ is the sum of all descents of $\pi$. 

An index pair $1 \le
i < j \le n$ forms an {\em inversion} of $\pi$ if $\pi(i) >
\pi(j)$. Let $\inv(\pi)$ be the number of inversions of $\pi$.

In order to study these statistics all at once, it is 
convenient to define the quantity
\[
I_\lambda(i,j) = \#\{\pi \in C_\lambda\mid \pi(i) > \pi(j)\}
\]
for $\lambda \vdash n$ and $1\leq i < j \leq n$.
\begin{lemma}\label{le:inv}
Suppose $\lambda\vdash n$. Let $p$ and $q$ denote the number of
$1$-parts and the number of $2$-parts, respectively, in
$\lambda$. Then,
\[
I_\lambda(i,j) =
\frac{\#C_\lambda}{2}\left(1 + \frac{2q-p(p-1)}{n(n-1)} +
  \frac{2(j-i-1)((n-p)(1-p) - 2q)}{n(n-1)(n-2)}\right).
\]
\begin{proof}
Fix $\lambda \vdash n$ and indices $1\leq i < j \leq n$. Consider the following subsets of $C_\lambda$:
\[
\begin{split}
T_1 &= \{\pi \in C_\lambda \mid \pi(i) = i\text{ and }\pi(j) = j\},\\
T_2 &= \{\pi \in C_\lambda \mid \pi(i) = j\text{ and }\pi(j) = i\},\\
T_3 &= \{\pi \in C_\lambda \mid \pi(i) = i\text{ and }\pi(j) = k\text{ for
  some }i<k<j\},\\
T_4 &= \{\pi \in C_\lambda \mid \pi(j) = j\text{ and }\pi(i) = k\text{ for
  some }i<k<j\},\\
T_5 &= \{\pi \in C_\lambda \mid \pi(i) = j\text{ and }\pi(j) = k\text{ for
  some }i<k<j\},\\
T_6 &= \{\pi \in C_\lambda \mid \pi(j) = i\text{ and }\pi(i) = k\text{ for
  some }i<k<j\},\\
T_7 &= C_\lambda \setminus(T_1\cup T_2\cup T_3\cup T_4\cup T_5\cup T_6).
\end{split}
\]
Thus, $C_\lambda$ is the disjoint union $C_\lambda = T_1 \uplus \cdots \uplus
T_7$.

Let $f_{i,j}:C_\lambda \to C_\lambda$ be the involution $\pi \mapsto
(i\, j)\pi(i\, j)$, where $(i\, j)$ denotes the transposition which
interchanges $i$ and $j$. Then, $f_{i,j}$ restricts to an involution
$T_7\to T_7$. This restriction has no fixed points since $T_1\cup T_2$
is the fixed point set of $f_{i,j}$. Moreover, for a permutation
$\pi\in T_7$, $(i, j)$ is an inversion if and only if it is not an inversion of
$f_{i,j}(\pi)$. 

Observing that $(i,j)$ is an inversion for all $\pi\in T_2\cup T_5\cup
T_6$, whereas it is a non-inversion for all $\pi\in T_1\cup T_3\cup
T_4$, we thus obtain
\[
\begin{split}
I_\lambda(i,j) &= \#T_2 + \#T_5 + \#T_6 + \frac{\#T_7}{2} \\
&= \frac{\#C_\lambda -\#T_1 + \#T_2 -\#T_3 - \#T_4 + \#T_5 + \#T_6}{2}.
\end{split}
\]
Computing
\[
\begin{split}
\#T_1 &= \frac{p(p-1)\#C_\lambda}{n(n-1)},\\
\#T_2 &= \frac{2q\#C_\lambda}{n(n-1)},\\
\#T_3 &= \#T_4 = \frac{p(n-p)(j-i-1)\#C_\lambda}{n(n-1)(n-2)},\\
\#T_5 &= \#T_6 = \frac{(n-p-2q)(j-i-1)\#C_\lambda}{n(n-1)(n-2)}\\
\end{split}
\]
yields the desired result.
\end{proof}
\end{lemma}

Next, we exploit the fact that several familiar permutation statistics
are obtained by taking appropriate sums of $I_\lambda(i,j)$. 

\begin{theorem}\label{th:inv}
The mean statistics associated with $\des$, $\maj$ and $\inv$ can
be written as the following linear combinations of irreducible characters:
\[
\begin{split}
\av{\des} &= \frac{n-1}{2}\chi^{(n)} - \frac{1}{n}\chi^{(n-1,1)} - \frac{1}{n}\chi^{(n-2,1,1)},\\
\av{\maj} &= \frac{n}{2}\av{\des} = \frac{n(n-1)}{4}\chi^{(n)} - \frac{1}{2}\chi^{(n-1,1)} - \frac{1}{2}\chi^{(n-2,1,1)},\\
\av{\inv} &= \frac{n(n-1)}{4}\chi^{(n)} - \frac{n+1}{6}\chi^{(n-1,1)} -
\frac{1}{6}\chi^{(n-2,1,1)}.\\
\end{split}
\]
\begin{proof}
Let $\lambda \vdash n$. Applying Lemma \ref{le:inv}, we obtain the identities
\[
\begin{split}
\av{\des}(\lambda) &=
\sum_{i=1}^{n-1}\frac{I_\lambda(i,i+1)}{\#C_\lambda} = \frac{n-1}{2}
+ \frac{q}{n} - \frac{1}{n}\bin{p}{2},\\
\av{\maj}(\lambda) &=
\sum_{i=1}^{n-1}\frac{i I_\lambda(i,i+1)}{\#C_\lambda} = \frac{n(n-1)}{4}
+ \frac{q}{2} - \frac{1}{2}\bin{p}{2},\\
\av{\inv}(\lambda) &=
\sum_{1\le i<j\le n}\frac{I_\lambda(i,j)}{\#C_\lambda} = \frac{n(n-1)}{4}
- \frac{p(p-1)}{12} - \frac{n(p-1)}{6} + \frac{q}{6},\\
\end{split}
\]
where $p$ and $q$ are as in Lemma \ref{le:inv}. That these equations
are equivalent to the asserted ones is readily shown using Lemma \ref{le:explicit}.
\end{proof}
\end{theorem}

\section{Some examples}\label{se:examples}
If $\Gamma$ is a union of conjugacy classes, we may combine Theorem
\ref{th:inner} with Theorem \ref{th:blambda} in order to explicitly
compute $\E_\Gamma(\stat, t)$ for any of the permutation statistics
$\stat$ which were studied in the previous sections. We conclude the
paper with a few sample computations of this kind.

\begin{example}\label{ex:inv}
Suppose $\Gamma = T \subseteq \Sn$ is the set of transpositions. The
contents $c_\lambda$ considered in Lemma \ref{le:content} make an appearance
as Theorem \ref{th:blambda} specializes to
\[
b_\lambda^{(t)} =
\frac{1}{n!}\left(\frac{\bin{n}{2}\chi^\lambda((2,
    1^{n-2}))}{f^\lambda}\right)^tf^\lambda = \frac{c_\lambda^tf^\lambda}{n!}.
\]
Combining this with the decompositions found in Theorem \ref{th:exc}
and Theorem \ref{th:inv}, we may invoke
Theorem \ref{th:inner} and show that a product of $t$ random transpositions has the following expected values of the
(weak) excedance number, descent number, major index and inversion number, respectively:
\[
\begin{split}
\E_T(\exc, t) &= \frac{n-1}{2}\left(1 - \left(1 - \frac{2}{n-1}\right)^t\right),\\
\E_T(\wexc, t) &= \frac{n+1}{2}\left(1 + \frac{n-1}{n+1}\left(1 - \frac{2}{n-1}\right)^t\right),\\
\E_T(\des, t) &= \frac{n-1}{2}\left(1 -
  \frac{2}{n}\left(1-\frac{2}{n-1}\right)^t -
  \frac{n-2}{n}\left(1-\frac{4}{n-1}\right)^t\right),\\
\E_T(\maj, t) &= \frac{n(n-1)}{4}\left(1 -
  \frac{2}{n}\left(1-\frac{2}{n-1}\right)^t - \frac{n-2}{n}\left(1-\frac{4}{n-1}\right)^t\right),\\
\E_T(\inv, t) &= \frac{n(n-1)}{4}\left(1 -
\frac{2(n+1)}{3n}\left(1-\frac{2}{n-1}\right)^t - \frac{n-2}{3n}\left(1-\frac{4}{n-1}\right)^t\right).\\
\end{split}
\]
\end{example}

The formula for $\E_T(\inv, t)$ obtained in Example \ref{ex:inv} was
previously found by Sj\"ostrand \cite[Theorem 5.1]{sjostrand}. 

We have refrained from
stating the explicit formula for $\E_T(\cyc_k,t)$ which follows from
Theorem \ref{th:cycles}. With $k=1$, it recovers results on fixed
points from \cite{jonsson}, whereas the sum over all $k$ leads to
expected cycle numbers that were computed in \cite{EH}. Also, note that
$\E_T(\cyc_n,t)/n$ is nothing but the probability that a product of $t$ random
transpositions forms an $n$-cycle. From that probability, one
easily derives the number of factorisations of an $n$-cycle
into $t$ transpositions. Working out the details, one recovers the
formula counting such factorisations which appears in Jackson \cite{jackson}.

\begin{example}
Let $\Gamma = C_{(n)}\subseteq \Sn$ be the set of $n$-cycles. Assuming now
that $t > 0$, Theorem \ref{th:blambda} and Lemma \ref{le:hook} show that $b_\lambda^{(t)} =
0$ unless we have a hook shape $\lambda = (n-j, 1^j)$. Moreover,
\[
b_{(n-j, 1^j)}^{(t)} = \frac{1}{n!}\left(\frac{(n-1)!(-1)^j}{\bin{n-1}{j}}\right)^t\bin{n-1}{j}.
\]
Suppose $k < n$ and consider the statistic $\cyc_k$ which counts
elements that belong to $k$-cycles. According to Theorem \ref{th:cycles},
exactly two terms (namely $\chi^{(n)}$ and
$(-1)^{k-1}\chi^{(n-k, 1^k)}$) which correspond to hook shapes appear in the expansion of
$\cyc_k$. Applying Theorem \ref{th:inner}, we
deduce that the expected number of $k$-cycles in a product of $t>0$
random $n$-cycles in $\Sn$ is
\[
\frac{1}{k}\E_{C_{(n)}}(\cyc_k,t) = \frac{1}{k} + \frac{(-1)^{k(t+1)-1}}{k\bin{n-1}{k}^{t-1}}.
\]
\end{example}

Let us, finally, consider an example where $\Gamma$ does not consist of a
single conjugacy class.

\begin{example}
Suppose $\Gamma$ is any union of conjugacy classes in which every
permutation has exactly one fixed point. In this case, Theorem
\ref{th:blambda} in conjunction with Lemma \ref{le:explicit} shows
that $b_{(n-1,1)}^{(t)} = 0$ for all $t > 0$. Thus, by Theorem
\ref{th:inner} and Theorem \ref{th:exc}, the expected number of
excedances of a product of $t>0$ elements from $\Gamma$ is
\[
\E_\Gamma(\exc, t) = \frac{n-1}{2}, 
\]
whereas for weak excedances one obtains
\[
\E_\Gamma(\wexc, t) = \frac{n+1}{2},
\]
independently of $t$.
\end{example}

\end{document}